\documentclass[12pt]{amsart}
\usepackage{amssymb,amsthm,amsmath,xypic,amscd}
\newtheorem{thm}{Theorem}[section]

\newtheorem{lem}[thm]{Lemma}
\newtheorem{prop}[thm]{Proposition}
\theoremstyle{definition}
\newtheorem{defn}[thm]{Definition}
\DeclareMathOperator{\supp}{\textsf{supp}\:}
\DeclareMathOperator{\Char}{char}

 \DeclareMathOperator{\Inf}{Inf}
 \DeclareMathOperator{\Res}{Res}

\newcommand{\csa}{$C^{*}$-algebra }
\newcommand{\csas}{$C^{*}$-algebras }

\newcommand{\norm}[1]{\left\Vert#1\right\Vert}
\newcommand{\abs}[1]{\left\vert#1\right\vert}

\newcommand{\Z}{\mathbb{Z}}
\newcommand{\C}{\mathbb{C}}
\newcommand{\Q}{\mathbb{Q}}
\newcommand{\F}{\mathbb{F}}
\newcommand{\K}{\mathcal{K}}
\newcommand{\A}{\mathbb{A}}
\newcommand{\R}{\mathbb{R}}

\newcommand{\E}{\underline{\text{E}}}
\newcommand{\Pc}{{\mathcal{P}_{c}}}

\newcommand{\Kt}[1]{K_{*}^{\textrm{top}}({#1})}

\begin{document}
\title{Local-global principle for the Baum-Connes conjecture with coefficients}
\author{Paul Baum, Stephen Millington and Roger Plymen}
\date{version 6: \today}
\begin{abstract} We establish the Hasse principle (local-global
principle) in the context of the Baum-Connes conjecture with
coefficients. We illustrate this principle with the discrete group
$GL(2,F)$ where $F$ is any global field.
\end{abstract}
\maketitle
\section{Introduction}
Let $G$ be a second countable locally compact Hausdorff
topological group. We shall say that $G$ satisfies $BCC$, or $BCC$
is true for $G$, if the Baum-Connes conjecture with coefficients
in an arbitrary $G-C^*$-algebra is true for $G$.

Our first result is the following permanence property:

\begin{thm}
\label{first} Let $G$ be the ascending union of open subgroups
$G_n$, and let $A$ be a $G-C^*$-algebra.  If each subgroup $G_n$
satisfies the Baum-Connes conjecture with coefficients $A$, then
$G$ satisfies the Baum-Connes conjecture with coefficients $A$.
\end{thm}

Our main application is to the following new permanence property.

\begin{thm}
\label{main} Let $F$ be a global field, $\A$ its ring of adeles,
$G$ a linear algebraic group defined over $F$.  Let $F_v$ denote a
place of $F$. If BCC is true for each local group $G(F_v)$ then
BCC is true for the adelic group $G(\A)$.
\end{thm}

Another application of Theorem 1.1 is the proof of the Baum-Connes
conjecture for reductive adelic groups \cite{BMP}.

To derive Theorem \ref{main} from Theorem \ref{first} we note that
the adelic group $G(\A)$ admits an ascending union of open
subgroups. We then make use of a crucial permanence property due
to Chabert-Echterhoff \cite{CE}, namely that BCC is stable under
direct product of finitely many groups.

If $G$ satisfies $BCC$, then any closed subgroup of $G$ also
satisfies $BCC$ \cite[Theorem 2.5]{CE}. Since $G(F)$ is a discrete
subgroup of $G(\mathbb{A})$, we have the following result:

\begin{thm}
\label{discrete}  If BCC is true for each local group $G(F_v)$
then BCC is true for the discrete group $G(F)$.
\end{thm}

There is, at present, a limited supply of local groups for which
BCC is known to be true.   Nevertheless, some examples are known.

For the group $SO(n,1)$, Kasparov \cite{K} proved that $\gamma =
1_G$ in the Kasparov representation ring $R(G) = KK_G(\C, \C)$.
For the group $SU(n,1)$, Julg-Kasparov \cite{JK} and
Higson-Kasparov \cite{HK} proved that $\gamma = 1_G$ in the ring
$R(G)$. For the group $Sp(n,1)$, Julg \cite{JJ} has recently
proved that, for any $G-C^*$-algebra $A$, the image of $\gamma$
via the map:
\[R(G) \to KK_G(A,A) \to KK(A \rtimes_{r}G, A \rtimes_{r} G) \to End \,K_*(A \rtimes_{r} G)\]
is the identity element even though $\gamma \neq 1$ in the ring
$R(G)$. Therefore the following rank-one Lie groups satisfy $BCC$:
\[SO(n,1),\; SU(n,1),\; Sp(n,1).\]

Let $F$ be a global field and consider the group $GL(2,F)$. At
each place $v$ of $F$ we have the local group $GL(2,F_v)$. The Lie
group $SL(2,\R)$ acts properly and isometrically on the
Poincar\'{e} disc, and so has the Haagerup property
\cite[p.159]{J}. The Lie group $SL(2,\C)$ acts properly and
isometrically on the Poincar\'{e} ball \cite[p.11]{E}, and so has
the Haagerup property \cite[p.159]{J}. If $v$ is a finite place
then the totally disconnected group $SL(2,F_v)$ acts properly and
isometrically on its tree \cite{S}, admits a proper $1$-cocycle,
and so has the Haagerup property \cite[p.155]{J}.

For each place $v$, the determinant map creates a short exact
sequence
\[
1 \to SL(2,F_v) \to GL(2,F_v) \to F_v^{\times} \to 1
\]
in which $SL(2,F_v)$ is a closed normal subgroup with the Haagerup
property, and the locally compact abelian group $F_v^{\times}$
certainly satisfies $BCC$. Therefore, by \cite[Corollary 3.14]{CE}
each local group $GL(2,F_v)$ satisfies $BCC$.

Let $\A_F$ be the adele ring attached to $F$. By Theorem
\ref{main}, we have

\begin{thm}
The adelic group $GL(2,\A_F)$ satisfies $BCC$.
\end{thm}

By Theorem \ref{discrete}, we have

\begin{thm}  The discrete group $GL(2,F)$, and each of its
subgroups, satisfies BCC.
\end{thm}

Our method is therefore an example of the Hasse principle
(local-global principle) applied to the {\it local} groups
$G(F_v)$ and the {\it discrete} group $G(F)$.

It is worth noting that if the Baum-Connes conjecture fails for
the discrete group $SL(n, \Z)$, then BCC fails for $SL(n, \R)$.

The discrete group $SL(2, F)$, with $F$ a global field, is
discussed in the Bourbaki seminar by Pierre Julg \cite[p. 160]{J}.
He proves that $SL(2,F)$ has the Haagerup property. It then
follows from a theorem of Higson-Kasparov \cite{HK} that $SL(2,
F)$ satisfies the Baum-Connes conjecture, see \cite[p.152]{J}.

In the course of our work, we find it necessary to use a model
$P_c(G)$ of the universal example for proper actions of $G$ which
is itself a direct limit of compact spaces. This model $P_c(G)$ is
paracompact, Hausdorff and separable, but not metrizable, and so
falls outside the discussion of proper actions in \cite{BCH}. We
have therefore to choose a different starting point for the theory
of proper actions: we use the definition of Bourbaki
\cite{bourbaki}.

This paper is a sequel to our Note \cite{BMP}. We thank Siegfried
Echterhoff, Pierre Julg, Ryszard Nest and Richard Sharp for
valuable conversations.

\section{Adelic groups}

 A \emph{local field} is
a non discrete locally compact topological field. It is shown in
\cite{We} that a local field $F$ must be of the following form. If
$\Char(F)=0$, then $F=\R$, $\C$ or a finite extension of $\Q_{p}$
for some prime $p$.   If $\Char(F)=p>0$, then $F$ is the field
$\mathbb{F}_{q}((X))$ of formal Laurent series (with finite tail)
in one variable with coefficients in a finite field
$\mathbb{F}_{q}$.

The fields $\R, \C$ are known as \emph{archimedean fields}. All
other local fields are known as \emph{nonarchimedean fields}. The
topology on a nonarchimedean field is always totally disconnected.

Let $\F_p(t)$ denote the field of fractions of the polynomial ring
$\F_p[t]$. A \emph{global field} is a finite extension of $\Q$, or
a finite extension of the function field $\F_p(t)$. A
\emph{completion} $(v, F_{v})$ of $F$ is a dense isomorphic
embedding $v$ of $F$ into a local field $F_{v}$. Two completions
$(v, F_{v})$, $(u, F_{u})$ are said to be equivalent if there is
an isomorphism $\rho$ of $F_{v}$ onto $F_{u}$ such that
$u=\rho\circ v$. A \emph{place} of $F$ is an equivalence class of
completions. We say the place $(v, F_{v})$ is \emph{infinite} if
$F_{v}$ is an archimedean field and \emph{finite} otherwise. If
$\Char(F)=p>0$, then $F$ has countably many finite places and
\emph{no} infinite places.   If $\Char(F)=0$ then $F$ has
countably many finite places and finitely many (but at least one)
infinite places.

\label{seclim} Suppose we have an ascending sequence of
topological spaces
\[X_{1}\subset X_{2}\subset X_{3}\subset\dots\]
Then we can give the union $X=\cup X_{n}$ the \emph{direct limit
topology}: that is a set is open in $X$ if and only if it has open
intersection with each $X_{n}$.

The following proposition admits a standard proof in point-set
topology.

\begin{lem}
\label{keyprop}

Let $X_{1}\subset X_{2}\subset X_{3}\subset\dots$ be a sequence of
$T_{1}$ (=points are closed) topological spaces and give
$X=\cup_{n=1}^{\infty}X_{n}$ the direct limit topology. Then any
compact subset of $X$ lies entirely within some $X_{n}$.
\end{lem}

Suppose that for each $n=1, 2, 3, \dots$ we have a locally
compact, second countable and Hausdorff topological group $G_{n}$,
such that $G_{m}$ is an \emph{open} subgroup of $G_{n}$ for $m\leq
n$:
\[G_{1}\subset G_{2}\subset G_{3}\subset\dots\]
Let $G=\cup_{n=1}^{\infty} G_{n}$ and furnish this with the direct
limit topology. Then $G$ is a topological group which is locally
compact, second countable and Hausdorff.

Any nonarchimedean local field $F_{v}$ contains a unique maximal
compact open subring $\mathcal{O}_{v}$. Let $S$ denote any finite
set of places of $F$ which contains all the infinite places.   By
an {\it adele} we mean an element $a = (a)_v$ of the product
$\prod_vF_v$ such that $a \in A_S = \prod_{v\in S} F_v \times
\prod_{v \notin S} \mathcal{O}_v$ for some $S$.   The adeles of
$F$ form a ring $\A_F$, addition and multiplication being defined
componentwise.   Each $A_S$ has its natural topology and $\A_F =
\cup_S A_S$ is topologized as the inductive limit with respect to
$S$.   There is an obvious embedding of $F$ in $\A_F$, by means of
which we identify $F$ with a subring of $\A_F$. The field $F$ is a
discrete cocompact subfield of the non-discrete locally compact
semisimple commutative ring $\A_F$.

Now suppose $G$ is a linear algebraic group defined over $F$. We
shall be interested in the \emph{adelic group} $G(\A_{F})$ of
$\A_{F}$-rational points of $G$. For a finite place $v$ of $F$ let
$G(\mathcal{O}_{v})$ denote the group $G(F_{v})\cap
GL_{n}(\mathcal{O}_{v})$. We set \[G_S = \prod_{v \in S} G(F_v)
\times \prod_{v \notin S} G(\mathcal{O}_v)\]

Then $G(\A)$ is equal, by definition, to the direct limit of the
groups $G_S$.   We now equip $G(\A)$ with the direct limit
topology, following Weil \cite[p.2]{WW}.   Then $G(\A)$ is a
locally compact second countable Hausdorff group.   The map $x
\mapsto (x,x,x,\ldots)$ embeds $G(F)$ as a {\it discrete} subgroup
of $G(\A)$.

\section{Proper actions and universal examples}
We recall that a topological space $X$ is {\it completely regular}
if it satisfies the following separation axiom: If $B$ is a closed
subset of $X$ and $p \in X \setminus B$ then there exists a
continuous function $f: X \to [0,1]$ such that $f(p) = 0$ and
$f(B) = \{1\}$. Let $G$ be a locally compact, Hausdorff, second
countable group. A \emph{$G$-space} is a topological space $X$
with a given continuous action of $G$ such that
\begin{itemize}
\item $X$ is completely regular and any compact subset of $X$ is metrizable.
\item The quotient space $X/G$ is paracompact and Hausdorff.
\end{itemize}

Note that the definition of a $G$-space in \cite{BCH} uses the
slightly more restrictive conditions that $X$ and $X/G$ be
metrizable. Nothing is altered in \cite{BCH} if this is replaced
throughout with the above relaxed conditions.

The following definition may be found in Bourbaki \cite[Definition
1, p.250, Prop.7, p.255]{bourbaki}.
\begin{defn}
\label{proper} The action of $G$ on a $G$-space $X$ is
\emph{proper} if given any two points $x, y\in X$ there are open
neighbourhoods $U_{x}, U_{y}$ of $x$ and $y$ respectively such
that the set
\[\{g\in G:gU_{x}\cap U_{y}\neq\emptyset\}\] has compact closure
in $G$. A proper $G$-space $X$ is said to be \emph{$G$-compact} if
the quotient $X/G$ is compact.
\end{defn}

In appendix A we give the full statement of Theorem 3.8 in Biller
\cite{biller}. This theorem, on the existence of slices,
reconciles the Bourbaki definition of proper actions with the
definition in \cite{BCH}.  In particular, if X is a proper G-space
such that X and and X/G are metrizable then X satisfies the
condition which was taken in [1] to be the definition of proper.
The reverse implication is also valid. If X is a metrizable
G-space which is proper in the sense of [1], then X is a proper
G-space.

Note that if $X$ is locally compact then this is equivalent to the
following condition: if $K_{1}, K_{2}$ are compact subsets of $X$
then the set
\[\{g\in G: gK_{1}\cap K_{2}\neq\emptyset\}\text{ is compact. }\]

It is easy to see any locally compact Hausdorff group acts
properly on itself. Also note if $X$ is a proper $G$-space and $Y$
is any $G$-invariant subset of $X$ then $Y$ is a proper $G$-space
when equipped with the subspace topology and obvious action of
$G$. We shall want to define $KK$ groups for the algebra
$C_{0}(X)$ of a $G$-compact proper $G$-space and the following
ensures that these algebras are separable.
\begin{prop}
\label{lcprop} Suppose $G$ acts properly on $X$ and $X/G$ is
compact. Then $X$ is locally compact and second countable.
\end{prop}
\begin{proof}
The space $X$ is locally compact by \cite[(e), p.310]{bourbaki}.
We show that $X=GS$ for some compact $S\subset X$.

For each $x\in X$ let $U_{x}$ be an open neighbourhood of $x$ with
compact closure. If $\pi:X\to X/G$ denotes the quotient map then
the collection $\pi(U_{x})$ cover $X/G$. We know $X/G$ is compact
and so get a finite subcover $\pi(U_{x_{1}}), \dots ,
\pi(U_{x_{n}})$. Now let $S=\cup_{i=1}^{n}\overline{U_{x_{i}}}$.
This is a compact set with $GS=X$.

By assumption any compact subset of $X$ is metrizable, and we may
deduce that $S$ is second countable when given the subspace
topology from $X$.

The action of $G$ on $X$ gives rise to a map
\[G\times S\to GS = X,\]
which is continuous, open and surjective. By assumption $G$ is
second countable and so clearly $X$ must be second countable.
\end{proof}

\begin{defn}
\label{ctofdefn} Let $X$ be a proper $G$-space. A \emph{cutoff
function} on $X$ is a function $c:X\to \R_{+}$ such that the
support of $c$ has compact intersection with $GK$ for any compact
subset $K$ of $X$ and
\[\int c(g^{-1}x) dg=1, \quad \text{ for each } x\in X\]
\end{defn}
Note the set of all such functions is convex.
\begin{prop}
\label{ctofexst} Let $X$ be a $G$-compact proper $G$-space. Then
there exists a cutoff function on $G$.
\end{prop}
\begin{proof} By Proposition \ref{lcprop} $X$ is locally compact
and $X=GS$ for some compact subset $S\subseteq X$. Take an open
neighbourhood $V$ of $S$ with compact closure and let $f$ be a
continuous function from $X$ into the interval $[0,1]$ which is
equal to $1$ on $S$ and zero on the complement of $V$, hence
compactly supported.
\\

Now for any $x\in X$ the function
\[f_{x}:G\to [0,1],\quad f_{x}(g)=f(g^{-1}x)\]
is continuous. Because $GS=X$ we may find a $g\in G$ with
$g^{-1}x\in S$. Recall $f=1$ on $S$ and so
$f_{x}(g)=f(g^{-1}x)=1$. As $f_{x}$ is continuous we can find a
neighbourhood of $g$ on which $f_{x}$ is non zero. Thus we may
conclude
\[0< \int f_{x}(g) dg\]
The action of $G$ on $X$ is proper so the set
\[\{g\in G:g\{x\}\cap\supp(f)\neq\emptyset\}\]
is compact for each $x\in X$. This implies that $f_{x}$ is
compactly supported and so
\[\int f_{x}(g) dg < \infty.\]
Thanks to the above we may define
\[c(x)=\frac{f(x)}{\int f_{x}(g) dg}\]
This is compactly supported and has the property that
$G\cdot\supp(c)=X$ and \[\int c(g^{-1}x) dg=1\] for each $x\in X$.
Clearly $c$ is a cutoff function on $X$.
\end{proof}
\begin{defn}
Let $X$, $Y$ be proper $G$-spaces. A continuous map $\varphi:X\to
Y$ is a called a \emph{G-map} if
\[\varphi(gx)=g\varphi(x)\text{ for all }g\in G, x\in X\]
Two $G$-maps are \emph{$G$-homotopic} if they are homotopic
through $G$-maps.
\end{defn}

\begin{defn}
\label{uedefn} A universal example for proper actions of $G$,
denoted $\E G$, is a proper $G$-space with the following property:
If $X$ is any proper $G$-space, then there exists a $G$-map
$f:X\to\E G$, and any two such maps are $G$-homotopic.
\end{defn}

Let $K$ be any compact subset of $G$.
The set of all probability measures on $K$
--- denoted $\mathcal{P}(K)$ --- is a separable compact Hausdorff space in
the topology induced from the weak* topology on $C(K)$. Recall
$\mu_{i}\to\mu$ in the weak* topology if and only if
\[\int f d\mu_{i}\to\int f d\mu\quad\text{for any } f\in
C(K).\]


If $K_{1}\subseteq K_{2}$ are compact subgroups of $G$ then define
the following map
\[\iota:\mathcal{P}(K_{1})\to\mathcal{P}(K_{2}),\quad
(\iota\mu)(U)=\mu(U\cap K_{1}) \text{ for each Borel set }
U\subseteq K_{2}\]

This is clearly injective and continuous. So $\iota$ is an
injective map from a compact space to a Hausdorff space and hence
gives rise to a homeomorphism between $\mathcal{P}(K_{1})$ and its
image in $\mathcal{P}(K_{2})$ with the subspace topology. To
simplify notation we identify $\mathcal{P}(K_{1})$ with its image
in $\mathcal{P}(K_{2})$.

As $G$ is locally compact and second countable it is clearly
$\sigma$-compact and by \cite[p.94, Prop. 15 and Cor. 2]{bourbaki}
we may find compact subsets $K_{1}\subset K_{2}\subset
K_{3}\subset \dots$ with $G=\cup_{i=1}^{\infty}K_{i}$ such that
\emph{any} compact subset of $G$ is contained within some $K_{i}$.
Now define
\[\mathcal{P}_{c}(G)=\cup_{i=1}^{\infty}\mathcal{P}(K_{i}).\]
Clearly $\Pc(G)$ consists of all compactly supported probability
measures on $G$. We topologize $\mathcal{P}_{c}(G)$ by giving it
the direct limit topology with respect to the sequence
$\mathcal{P}(K_{i})$.  This topology is independent of the
sequence $K_{i}$.

The group $G$ acts on $\mathcal{P}_{c}(G)$ by setting
\[(g\mu)(U)=\mu(g^{-1}U)\;\text{ for any Borel set } U\subseteq G.\]

\begin{lem}
\label{ueexst}  The action of $G$ on $\mathcal{P}_{c}(G)$ is
continuous.
\end{lem}

\begin{proof}

Let $(g_{\alpha}, \mu_{\alpha})\to (g, \mu)$, $\alpha\in A$, be a
convergent net in $G\times \mathcal{P}_{c}(G)$. By definition of
the direct limit topology on $\mathcal{P}_{c}(G)$, we must have
$\mu_{\alpha}\to\mu$ in $\mathcal{P}(K_{i})$ for some $i$. Given
$f\in C(K_{i})$ and $\epsilon>0$ we can find $\alpha_{0}$ with the
property that
\[\alpha>\alpha_{0}\Rightarrow\abs{\int f\; d\mu_{\alpha}-\int f\;
d\mu}\leq \epsilon/2 \] Furthermore $g_{\alpha}\to g$ and so we
may choose $\alpha_{0}$ in such a way that we also have
\[\alpha>\alpha_{0}\Rightarrow\norm{f_{g_{\alpha}}-f_{g}}\leq \epsilon/2\]
where $f_{g}$ denotes the function $x\mapsto f(g^{-1}x)$. Finally
note
\begin{eqnarray*}
\abs{\int f\; d(g_{\alpha}\mu_{\alpha})-\int f\; d(g\mu)} & = &
\abs{\int f_{g_{\alpha}}\; d\mu_{\alpha}-\int f_{g}\; d\mu}\\ &
\leq & \abs{\int (f_{g_{\alpha}}-f_{g})\; d\mu_{\alpha}} +
\abs{\int f_{g}\; d\mu_{\alpha}-\int f_{g}\; d\mu}\\ & \leq &
\epsilon/2+\epsilon/2
\end{eqnarray*}
This shows that $g_{\alpha}\mu_{\alpha}\to g\mu$ in
$\mathcal{P}(K_{i})$ and therefore $g_{\alpha}\mu_{\alpha}\to
g\mu$ in $\mathcal{P}_{c}(G)$.
\end{proof}

\begin{lem}
The action of $G$ on $\Pc(G)$ is proper.
\end{lem}
\begin{proof}
Take $\mu\in\Pc(G)$ and let $f_{\mu}$ be a continuous compactly
supported function $f_{\mu}:G\to [0,1]$ with $f_{\mu}\equiv 1$ on
$\supp(\mu)$.

Then the set
\begin{eqnarray*}
U_{\mu} & = & \left\{\lambda\in \Pc(G):\abs{\int f_{\mu}\;
d\mu-\int f_{\mu}\; d\lambda}\leq 1/2\right\}\\ & = &
\left\{\lambda\in\Pc(G):\int f_{\mu}\; d\lambda>1/2\right\}
\end{eqnarray*}
is an open neighbourhood of $\mu$ in $\Pc(G)$.

Now take any $\mu, \nu\in\Pc(G)$ and assume $gU_{\mu}\cap
U_{\nu}\neq\emptyset$ for some $g\in G$. Indeed let $\lambda\in
gU_{\mu}\cap U_{\nu}$. Then
\[\int g^{-1}f_{\mu}\; d\lambda>1/2\quad\text{and}\quad\int f_{\nu}\;
d\lambda>1/2\] If $\supp(g^{-1}f_{\mu})$ and $\supp(f_{\nu})$ are
disjoint then we have $0\leq g^{-1}f_{\mu}+f_{\nu}\leq 1$. However
by the above
\[\int g^{-1}f_{\mu}+f_{\nu}\;d\lambda >1\]
clearly contradicting the fact that $\lambda$ is a probability
measure. Hence we may conclude that
$g\supp(\mu)\cap\supp(\nu)\neq\emptyset$ and so
\[\{g\in G:gU_{\mu}\cap U_{\nu}\neq\emptyset\} \subset
\{g\in G:g\supp(\mu)\cap\supp(\nu)\neq\emptyset\}.\] Now both
$\mu$ and $\nu$ are compactly supported and since any group acts
properly on itself the larger set here is compact.
\end{proof}
\begin{thm}  The space $\Pc(G)$ is a universal example for proper actions of $G$.
\end{thm}
\begin{proof}
Let $X$ be \emph{any} proper $G$-space, we aim to show there
exists a $G$-equivariant map $X\to\Pc(G)$. Take any $x\in X$. By
\cite[Theorem 3.8]{biller} there is a $G$-invariant open
neighbourhood $U_{x}$ of $x$, a compact subgroup $H$ of $G$, and a
$G$-equivariant map
\[\rho:U_{x}\to G/H.\]
Let $\mu_{H}$ denote Haar measure on the compact subgroup $H$,
normalized to have total mass $1$. There is an obvious
$G$-equivariant map
\[\Psi:G/H\to \Pc(G),\quad gH\mapsto g\cdot\mu_{H}\]
and let
\[\theta_{x}=\Psi\circ\rho:U_{x}\to\Pc(G).\]

Now $X$ may be covered by such neighbourhoods and if $\pi$ denotes
the quotient map $X\to X/G$ then $\{\pi(U_{x})\}$ is an open cover
of $X/G$. Recall that by definition $X/G$ is paracompact and
Hausdorff and so there is a locally finite partition of unity
subordinate to the cover $\{\pi(U_{x})\}$. Precisely there are
continuous maps
\[\omega_{x}:X/G\to [0,1],\quad\text{for each } x\in X\]
with $\supp(\omega_{x})\subseteq \pi(U_{x})$ and for each $[y]\in
X/G$ we have
\[\omega_{x}([y])=0 \text{ for almost all }x\in X\text{ and }
\sum_{x\in X} \omega_{x}([y])=1.\] Now
\[\Xi:X\to\Pc(G),\quad \Xi(y)=\sum_{x\in X}
\omega_{x}(\pi(y))\theta_{x}(y)\] is the required map.  Note that
(as remarked above) $P_c(G)$ is a convex set and this convexity is
being used in the construction of $\Xi$.

Finally, if $\varphi_{1}, \varphi_{2}:X\to\Pc(G)$ are
$G$-equivariant maps then they are $G$-homotopic via
\[t\varphi_{1}+(1-t)\varphi_{2}, \quad t\in[0,1].\]
So $\Pc(G)$ is a universal example for $G$.
\end{proof}

\label{dirac} \emph{Baum-Connes conjecture with coefficients}. If
$A$ is any $G$-\csa then we may define

\[K^{\text{top}}_*(G,A) = \varinjlim_{\substack{\text{G-invariant}\\ \text{G-compact}\\ Z\subseteq \E G
}}KK_{*}^{G}(C_{0}(Z), A)\]

We say a group $G$ satisfies the Baum-Connes conjecture with
coefficients if for every $G$-\csa $A$ the map
\[\mu_{A}:K_{*}^{\text{top}}(G, A) \to K_{*}(A\rtimes_{r} G)\]
is an isomorphism.

\section{$K$-theory for ascending unions of groups}

\begin{thm}
\label{Tmap} Let $H$ be an open subgroup of $G$. Then the
inclusion of $H$ in $G$ determines a homomorphism of abelian
groups
\[\mathcal{T}_{H}^{G}:K_{*}(A\rtimes_{r} H)\to K_{*}(A\rtimes_{r} G).\]
Furthermore suppose $G_{1}\subset G_{2}\subset G_{3}\subset\dots$
is an ascending sequence of open subgroups , then there is an
inductive system of abelian groups
\[
\begin{CD}
K_{*}(A\rtimes_{r} G_{1}) @>\mathcal{T}_{G_{1}}^{G_{2}}>>
K_{*}(A\rtimes_{r}G_{2}) @>\mathcal{T}_{G_{2}}^{G_{3}}>>
K_{*}(A\rtimes_{r}G_{3}) @>\mathcal{T}_{G_{3}}^{G_{4}}>> \dots
\end{CD}
\] If $G=\cup G_{n}$ we have
\[K_{*}(A\rtimes_{r} G)=\varinjlim_{n} K_{*}(A\rtimes_{r} G_{n}).\]
\end{thm}
\begin{proof}
Any open subgroup $H$ of $G$ is also closed.  Let $H$ be an open
subgroup of $G$ and fix a Haar measure $dg$ on $G$. Since $H$ is
an open subgroup of $G$, the restriction of $dg$ to $H$ is a Haar
measure on $H$. To simplify notation we write $dg$ to signify this
choice of Haar measure on both $G$ and $H$.
\\
\\
For $f\in C_{c}(H,A)$ define

\[\iota f:G\rightarrow A, \quad (\iota f)(x)=\left\{
\begin{array}{ll}
f(x) & \text{if $x\in H$}\\
0 & \text{otherwise}.
\end{array}
\right.
\]
Then $\iota$ defines a map from $C_{c}(H,A)$ into $C_{c}(G,A)$ as
$H$ is clopen in $G$. $\iota$ is clearly additive and furthermore
for $f, g\in C_{c}(H,A)$
\begin{eqnarray*}
\iota f*\iota g(x) & = & \int_{G}\iota f(y) \alpha_y(\iota
g(y^{-1}x)) dy\\ & = & \int_{H}f(y)\alpha_y(\iota g(y^{-1}x)) dy\\
& = & \left\{
\begin{array}{ll}
\int_{H}f(y)\alpha_y(g(y^{-1}x)) dy & \text{if $x\in H$}\\ 0 &
\text{otherwise}.
\end{array}
\right. \\
& = & \iota(f*g)(x).
\end{eqnarray*}
Also for $f\in C_{c}(H,A)$
\begin{eqnarray*}
(\iota f)^{*}(x) & = & \left\{
\begin{array}{ll}
\Delta_{G}(x^{-1})\alpha_x(f(x^{-1})^*) & \text{if $x\in H$}\\ 0 &
\text{otherwise}
\end{array}
\right. = (\iota f^{*})(x)
\end{eqnarray*}
This follows from the observation that
$\Delta_{H}(x)=\Delta_{G}(x)$ for $x\in H$.

We now follow \cite[II.C, p.172]{C}.  Let $\norm{~}_{A \rtimes_{r}
G}$ signify the norm obtained on $C_{c}(G,A)$ from the left
regular representation on the $C^*$-module

\[E = L^{2}(G, dg) \otimes A.\]

We shall also need the $C^*$-module

\[F = L^2(H, dg) \otimes A.\]

If $f\in C_{c}(H,A)$ then
\begin{eqnarray*}
{\norm{\iota f}}_{A \rtimes_{r} G} & =& \sup\left\{{\norm{\iota
f*\xi}}_{E}:\xi \in E, \norm{\xi}_{E}\leq 1\right\}\\ & = &
\sup\left\{{\norm{f*\xi}}_{F}:\xi \in F, \norm{\xi}_{F}\leq
1\right\}
\end{eqnarray*}
So $\iota$ is an isometric $^{*}$-algebra map from $C_{c}(H,A)$ to
$C_{c}(G,A)$. We can then complete this map to obtain an injective
\csa morphism
\[\iota: A\rtimes_{r} H \to A\rtimes_{r}G.\]
The map $\mathcal{T}_{H}^{G}$ is then defined by functoriality:
\[\mathcal{T}_{H}^{G}=\iota_{*}:K_{*}(A\rtimes_{r}H)\to K_{*}(A\rtimes_{r}G).\]
Now let $G_{1}\subset G_{2}\subset G_{3}\subset\dots$ be an open
ascending sequence of groups and let $G=\cup G_{n}$. Fixing a Haar
measure on $G$ prescribes a Haar measure on each $G_{n}$. For any
$m<n$, $G_{m}$ is an open subgroup of $G_{n}$ and so by the above
we may define a map \[\iota_{m}^{n}: A\rtimes_{r} G_{m}\to
A\rtimes_{r} G_{n}.\] Indeed we may form an inductive system of
\csas $(A\rtimes G_{n}, \iota_{m}^{n})$,
\[
\begin{CD}
A \rtimes_{r} G_{1} @>\iota_{1}^{2}>> A \rtimes_{r} G_{2}
@>\iota_{2}^{3}>> A \rtimes_{r}G_{3} @>\iota_{3}^{4}>> \dots
\end{CD}
\]
We aim to show that the $C^*$-inductive limit $\varinjlim_{n}
A\rtimes_{r} G_{n}$ of this inductive system is equal to the
algebra $A\rtimes_{r} G$. For any $n$, $G_{n}$ is an open subgroup
of $G$ and let $\iota_{n}: A\rtimes_{r} G_{n}\rightarrow
A\rtimes_{r} G$ be defined as above. Since each $\iota_{n}$ is
injective, all that remains to be demonstrated is that
\[\bigcup_{n}\iota_{n}(A \rtimes_{r} G_{n})\] is dense in $A \rtimes_{r} G$.

Taking $\varepsilon>0$ and $x\in A \rtimes_{r} G$ we can find
$f\in C_{c}(G,A)$ with  ${\norm{f-x}}_{A \rtimes_{r}
G}<\varepsilon$. The support of $f$ is compact and by Lemma
\ref{keyprop} must lie totally within some $G_{n}$. Hence $f\in
C_{c}(G_{n},A)$ for some $n$ and so $\cup_{n}\iota_{n}(A
\rtimes_{r} G_{n})$ is dense in $A \rtimes_{r} G$.   Hence $A
\rtimes_{r} G$ is the $C^*$-inductive limit of the $C^*$-algebras
$A \rtimes_{r} G_n$ and we therefore have

\[K_{*}(A \rtimes_{r} G)  =  \varinjlim_{n}K_{*}(A
\rtimes_{r} G_{n}).\]

\end{proof}

\section{Equivariant $K$-homology for ascending sequences of
groups}
\begin{thm} \label{Rmap} Let $H$ be an open subgroup of $G$. Then the inclusion of
$H$ in $G$ determines a homomorphism of abelian groups
\[\mathcal{R}_{H}^{G}:\Kt{H,A}\to \Kt{G,A}.\]
Furthermore suppose $G_{1}\subset G_{2}\subset G_{3}\subset\dots$
is an open ascending sequence of groups , then there is an
inductive system of abelian groups
\[
\begin{CD}
\Kt{G_{1},A} @>\mathcal{R}_{G_{1}}^{G_{2}}>> \Kt{G_{2},A}
@>\mathcal{R}_{G_{2}}^{G_{3}}>> \Kt{G_{3},A}
@>\mathcal{R}_{G_{3}}^{G_{4}}>> \dots \end{CD}
\] If $G=\cup G_{n}$ then we have
\[\Kt{G,A}=\varinjlim_{n} \Kt{G_{n},A}\]
\end{thm}

In the course of proving this result, we shall make use of the
{\it reciprocity isomorphism} \cite[p.157]{CE} in equivariant
$KK$-theory. Let $H$ be an open subgroup of $G$, let $A$ be an
$H-C^*$-algebra, let $Ind^G_H A$ denote the induced algebra and
let $B$ be a $G-C^*$-algebra. Then we have the reciprocity
isomorphism:

\[Inf^G_H : KK^H_*(A,B) \cong KK^G_*(Ind^G_H A, B).\]

If $A$ is commutative we have $A \cong C_0(X), Ind^G_H A \cong
C_0(G \times_{H}X)$ and so we have

\[Inf_X :K_{*}^{H}(C_{0}(X),B)\to K_{*}^{G}(C_{0}(G\times_{H}X), B).\]

\begin{lem}
\label{induced} Let $H$ be an open subgroup of $G$ and let $X$ be
any $H$-compact subset of $\Pc(H)\subset\Pc(G)$ then
$G\times_{H}X\cong G\cdot X$ as $G$-spaces.
\end{lem}
\begin{proof}
Recall the definition of the space $G\times_{H}X$. The group $H$
acts on the product $G\times X$ by setting $h\cdot(g, x
)=(gh^{-1}, hx)$ and $G\times_{H}X$ is the quotient. The action of
$G$ on $G\times_{H} X$ is given by $g'\cdot[g, x]=[g'g, x]$, where
$[g, x]$ is the equivalence class of the pair $(g, x)$.
\\
\\
The map $F_X$ is defined as follows:
\[
F_{X}:G\times_{H}X  \to  G\cdot X, \; [g, x] \mapsto gx.
\]

The map $F_{X}$ is clearly surjective and $G$-equivariant. To show
this map is injective take $[g_{1}, x_{1}], [g_{2}, x_{2}]\in
G\times_{H} X$ with $g_{1}x_{1}=g_{2}x_{2}$. Recall that $x_{1},
x_{2}$ are in fact probability measures on $H$ and so
\[1=x_{1}(H)=g_{1}^{-1}g_{2}x_{2}(H)=x_{2}(g_{2}^{-1}g_{1}H)\]
As $x_{2}$ is a measure in $\Pc(H)$ it will clearly be equal to
zero on the coset $g_{2}^{-1}g_{1}H$ unless $g_{2}^{-1}g_{1}\in
H$. Now $g_{2}^{-1}g_{1}\cdot (g_{1}, x_{1}) =
(g_{1}g_{1}^{-1}g_{2}, g_{2}^{-1}g_{1}x_{1}) = (g_{2}, x_{2})$ as
required.
\\\\
It now remains to show this map is a homeomorphism. Let $\pi$
denote the quotient map $G\times X\to G\times_{H} X$ and
$\theta:G\times X\to G\cdot X$ denote the map given by the action
of $G$ on $X$.  Both $\pi$ and $\theta$ are open, continuous maps.
Since $\theta = F_X \circ \pi$, this implies that $F_X$ is open
and continuous.
\end{proof}

Let $\mathcal{R}_{H, X}^{G}$ denote the composition
\[\begin{CD} KK_{*}^{H}(C_{0}(X), A) @>\Inf_{X}>>
KK^G_*(C_{0}(G\times_{H} X), A) @>(F_{X})_{*} >>
KK_{*}^{G}(C_{0}(G\cdot X), A)
\end{CD}\]

If $X$ and $Y$ are $H$-compact subsets of $\Pc(H)$ with $X\subset
Y$ then the following diagram commutes,
\[
\begin{CD}
KK_{*}^{H}(C_{0}(X), A) @> >> KK_{*}^{H}(C_{0}(Y), A)\\ @V
\mathcal{R}_{H, X}^{G} VV @V \mathcal{R}_{H, Y}^{G} VV\\
KK_{*}^{G}(C_{0}(G\cdot X), A) @>
>> KK_{*}^{G}(C_{0}(G\cdot Y), A)\\
\end{CD}
\]
where each horizontal map is given by inclusion of the spaces
involved.

Now $\Pc(H), \Pc(G)$ are universal examples for $H, G$
respectively. If $X$ is an $H$-compact subset of $\Pc(H)$ then
$G\cdot X$ is a $G$-compact subset of $\Pc(G)$.  Due to the above
the following map is well defined

\[\mathcal{R}_{H}^{G}:\Kt{H,A}\to\Kt{G,A}.\]

Now let $G_{1}\subset G_{2}\subset G_{3}\subset\dots$ be an
ascending sequence of open subgroups and let $G=\cup G_{n}$. There
is then an inductive system of abelian groups
\[
\begin{CD}
\Kt{G_{1},A} @>\mathcal{R}_{G_{1}}^{G_{2}}>> \Kt{G_{2},A}
@>\mathcal{R}_{G_{2}}^{G_{3}}>> \Kt{G_{3},A}
@>\mathcal{R}_{G_{3}}^{G_{4}}>> \dots \end{CD}
\]
and the maps
\[\mathcal{R}_{G_{n}}^{G}:\Kt{G_{n},A}\to\Kt{G,A}.\]

\begin{lem}
There exist compact sets $\Delta_1 \subset \Delta_2 \subset
\Delta_3 \subset \dots $ such that $\Delta_n \subset G_n$ and
$\cup_n \, Interior (\Delta_n) = G$.  Set $Z^{n,m} = G_n \cdot
P_c(G_n \cap \Delta_m)$.   Then

\begin{itemize}
\item $Z^{n,m} \subset P_c(G_n)$
\item $Z^{n,m}$ is preserved by $G_n$ and is $G_n$-compact
\item $Z^{n,m} \subset Z^{n,m+1}$
\item $\cup_m Z^{n,m} = P_c(G_n)$
\item $G_{n+1}\cdot Z^{n,m} \subset Z^{n+1,m}$
\end{itemize}

\end{lem}
\begin{proof}
Let $W_1, W_2, W_3, \dots$ be a countable basis for the topology
of $G$.   From this list $W_1, W_2, W_3, \dots$ delete all $W_j$
such that $\overline{W_j}$ is not compact.  Call the remaining
list $V_1, V_2, V_3, \dots $.   Let $U$ be any non-empty open set
in $G$ and choose $p \in U$.   Since $G$ is locally compact there
exists an open set $\Lambda$ in $Y$ with $p \in \Lambda$ and
$\overline{\Lambda}$ compact.   Consider $\Lambda \cap U$.   Now
$\Lambda \cap U$ is a union of $W_j$'s.   But any $W_j$ which is
contained in $\Lambda \cap U$ must in fact be a $V_k$ because
$\Lambda \cap U$ has compact closure.  $\Lambda \cap U$ is a union
of $V_k$'s (now let $p$ vary throughout $U$).

Now let $\Sigma_n = \overline{V_1 \cup V_2 \cup \dots \cup V_n}$
and $\Delta_n = G_n \cap \Sigma_n$.

Any compact set in $G$ is contained in some $\Delta_n$.

\end{proof}

\begin{lem}
\label{abelian} Let $\{A^{m,n}\}$ be a commutative diagram of
abelian groups in which the typical commutative square is

\[\begin{CD}
A^{n,m} @>>> A^{n,m+1}\\ @VVV @VVV\\ A^{n+1,m} @>>> A^{n+1, m+1}
\end{CD}\]
with $n,m = 1,2,3, \dots$.  Then there is a canonical isomorphism
of abelian groups:
\[
lim_{m \to \infty}(lim_{n \to \infty}\, A^{m,n}) \cong lim_{n \to
\infty}(lim_{m \to \infty}\, A^{m,n}).
\]
\end{lem}

\begin{proof}
Each side is canonically isomorphic to the direct limit
$\varinjlim A^{n,n}$ of the directed system $\{A^{n,n}\}$.
\end{proof}

\begin{thm}
$lim_{n \to \infty} \Kt{G_n,A} = \Kt{G,A}.$
\end{thm}

\begin{proof}
Consider the commutative diagram of abelian groups in which the
typical commutative square is:

\[\begin{CD}
\K_j^{G_n}(Z^{n,m},A) @>>> K_j^{G_n}(Z^{n,m+1},A)\\ @V \rho_n VV
@V \rho_n VV\\ K_j^{G_{n+1}}(Z^{n+1,m},A) @>>>
K_j^{G_{n+1}}(Z^{n+1,m+1},A)
\end{CD}\]

with $n,m = 1,2,3,\dots$.   Each horizontal arrow
\[
K^{G_n}_j(Z^{n,m},A) \to K_j^{G_n}(Z^{n,m+1},A)
\]

is the map of abelian groups determined by the inclusion $Z^{n,m}
\to Z^{n,m+1}$.   Each vertical map $\rho_n : K^{G_n}_j(Z^{n,m},A)
\to K^{G_{n+1}}_j(Z^{n+1,m},A)$ is the map
$\mathcal{R}^{G_{n+1}}_{G_n}$ followed by the map of abelian
groups induced by the inclusion $G_{n+1} \cdot Z^{n,m} \to
Z^{n+1,m}$.

We will write
\[
A^{n,m} = K^{G_n}_j(Z^{n,m}, A).
\]

Each $Z^{n,m}$ is $G_n$-compact, $\cup_m\,Z^{n,m} = P_c(G_n)$ and
any $G_n$-compact set in $P_c(G_n)$ is contained in some
$Z^{n,m}$.  Taking the direct limit along the $n$th row, we have

\[
lim_{m \to \infty} \,A^{m,n} = K^{top}_j(G_n,A).
\]

If we now take the direct limit in a vertically downward
direction, we obtain
\[lim_{n \to \infty} \,K^{top}_j(G_n,A).
\]



Now we fix attention on the $m$th column.   If $n \geq m$ then
$G_n \supset \Delta_m$ and so $G_n \cap \Delta_m = \Delta_m$. Then
we have
\[
A^{n,m} = K^{G_n}_j(G_n \cdot P_c(\Delta_m), A).
\]
Now the $G_n$-saturation of $P_c(\Delta_m)$ is equal to the
$G$-saturation of $P_c(\Delta_m)$.  Now we apply the reciprocity
isomorphism and we have
\[
A^{n,m} \cong K^G_j(G \cdot P_c(\Delta_m), A).
\]

if $n \geq m$.   Therefore the $m$th column {\it stabilizes} as
soon as $n \geq m$.   Therefore the direct limit down the $m$th
column is given by
\[
lim_{n \to \infty} A^{m,n} = K^G_j(G \cdot P_c(\Delta_m), A)).
\]
Now the sets $G \cdot P_c(\Delta_m)$ are cofinal in $G$-compact
sets in $P_c(G)$.   Taking the direct limit in the horizontal
direction, we have
\[
lim_{m \to \infty} lim_{n \to \infty} A^{m,n} \cong
K_j^{top}(G,A).
\]
By Lemma \ref{abelian} we have
\[
lim_{n \to \infty} K^{top}_j(G_n,A) \cong K^{top}_j(G,A).
\]
\end{proof}

\section{Adelic groups}
For $H$ an open subgroup of $G$ we have constructed a homomorphism
from $\Kt{H,A}$ to $\Kt{G,A}$, and likewise for the $K$ theory of
the reduced crossed product $C^{*}$-algebras. We wish to check
that these maps are compatible with the Baum--Connes $\mu$ map, i.
e. that the following diagram commutes.
\[
\xymatrix{
  \Kt{H,A} \ar[d]_{\mathcal{R}_{H}^{G}} \ar[r]^{\mu_{H}} & K_{*}(A \rtimes_{r} H)
  \ar[d]^{\mathcal{T}_{H}^{G}} \\
  \Kt{G,A} \ar[r]^{\mu_{G}} & K_{*}(A \rtimes_{r} G)}
\]
As a first step, we prove that the reciprocity isomorphism $Inf_X$
is compatible with the Baum-Connes $\mu$ map.

\begin{lem}
\label{*} Let $A$ be a $G-C^*$-algebra, let $H$ be an open
subgroup of $G$ and let $X$ be a locally compact proper
$H$-compact $H$-space. Then the following diagram commutes:
\[\begin{CD}
KK_{*}^{H}(C_{0}(X), A) @> \mu_{H} >> K_{*}(A \rtimes_{r}H)\\ @V
Inf_X VV @V \mathcal{T}_{H}^{G} VV\\
KK_{*}^{G}(C_{0}(G\times_{H}X), A) @> \mu_{G} >> K_{*}(A
\rtimes_{r}G)
\end{CD}\]
\end{lem}

\begin{proof}
The inverse of the reciprocity isomorphism $\Inf$ is the {\it
compression isomorphism} \cite[p. 157]{CE} and is given by the
composition $i_{*}\circ \text{Res}_{X}$, where $\text{Res}_{X}$ is
the obvious restriction map $KK_{*}^{G}(C_{0}(X), A)\to
KK_{*}^{H}(C_{0}(X), A)$ and $i$ is the inclusion
$C_{0}(X)\hookrightarrow C_{0}(G\times_{H} X)$ given by
\[i(f)[g, x]=\left\{
\begin{array}{ll}
f(gx) & \text{if $g\in H$}\\ 0 & \text{otherwise}.
\end{array}\right.
\]
Each of these maps is clearly functorial.

There is a commutative diagram
\[\begin{CD}
KK_{*}^{H}(C_{0}(X), A) @> j_{H} >> KK_{*}(C_{0}(X)\rtimes_{r} H,
A \rtimes_{r} H)\\ @A i_{*} AA @A i_{*}' AA\\
KK_{*}^{H}(C_{0}(G\times_{H}X), A) @> j_{H} >>
KK_{*}(C_{0}(G\times_{H}X)\rtimes_{r} H, A \rtimes_{r}H)\\ @A
\Res_{X} AA @V p_{*}\circ q_{*} VV\\
KK_{*}^{G}(C_{0}(G\times_{H}X), A) @> j_{G}
>> KK_{*}(C_{0}(G\times_{H}X)\rtimes_{r} G, A \rtimes_{r} G)
\end{CD}\]
in which $j_H, j_G$ are {\it descent} homomorphisms. Here $p$
denotes the map
\[p:C_{0}(G\times_{H}X)\rtimes_{r} G\to C_{0}(G\times_{H}X)\rtimes_{r} H\]
induced from the obvious restriction map $C_{c}(G, A)\to C_{c}(H,
A)$, with $A=C_{0}(G\times_{H}X)$. And $q$ denotes the map
\[q: A \rtimes_{r} H\to A \rtimes_{r} G\]
of Theorem \ref{Tmap}.

Now let $c$ be a cutoff function on the proper $H$-space $X$, we
claim $i(c)\in C_{0}(G\times_{H}X)$ is a cutoff function on the
proper $G$-space $G\times_{H}X$. To see this take any $[g_{0},
x]\in G\times_{H}X$. Then by definition we have
\[i(c)([g_{0}, x]) =  0\] unless
$g_{0}\in H$ in which case \[i(c)[g_{0}, x]=c(g_{0}x).\] Up to a
normalizing factor between the Haar measures on $H$ and $G$
\begin{eqnarray*} \int_{G}
i(c)(g^{-1}[g_{0}, x])\;dg & = & \int_{G}i(c)(g^{-1}g_{0}^{-1}[g_{0}, x])\;dg \\
& = & \int_{G} i(c)([g^{-1},x])\;dg  \\
& = & \int_{H} c(h^{-1}x)\;dh = 1.
\end{eqnarray*}
If $K$ is any compact subset of $G\times_{H}X$ and if $F$ denotes
the homeomorphism between $G\times_{H}X$ and $G\cdot X$ then
$F(GK)=G\cdot F(K)$ and $F(K)$ is compact. Also note
$F(\supp(i(c))\subseteq H\cdot\supp(c)$ and so
\[F(GK\cap\supp(i(c)))\subseteq G\cdot F(K)\cap
H\cdot\supp(c)\subseteq H\cdot F(K)\cap\supp(c)\] which is compact
and so $GK\cap\supp(i(c))$ is compact. So we have shown $i(c)$ is
a cutoff function on $G\times_{H}X$.

Let $\lambda_{X}$ denote the projection in the twisted convolution
algebra $C_{c}(H\times X)$ arising from the cutoff function $c$:

\[\lambda_Z(g,x) =
c_Z(x)^{1/2}c_Z(g^{-1}x)^{1/2}\Delta(g)^{-1/2}\]

and let $\lambda_{G\times_{H}X}$ denote the projection in
$C_{c}(G\times G\times_{H} X)$ arising from the cutoff function
$i(c)$.

Then $p(\lambda_{G\times_{H}X})$ is simply the restriction of
$\lambda_{G\times_{H}X}$ to $H\times G\times_{H} X$, and for any
$h$ in $G$ and any $[g, x]\in G\times_{H}X$
\begin{eqnarray*}
p(\lambda_{G\times_{H}X})(h, [g, x]) & = & i(c)[g,
x]^{1/2}i(c)(h^{-1}[g, x])^{1/2}\Delta_{G}(h)^{-1/2}\\\\
& = & \left\{\begin{array}{ll}
c(gx)^{1/2}c(h^{-1}gx)\Delta_{H}(h)^{-1/2} & \text{if $g\in H$}\\
0 & \text{otherwise}.
\end{array}\right.\\\\
& = & i'(\lambda_{X})
\end{eqnarray*}

So the following diagram commutes
\[\begin{CD}
KK_{*}(C_{0}(X)\rtimes_{r} H, A \rtimes_{r}H) @>
[\lambda_{X}]\otimes\cdot
>> KK_{*}(\C, A \rtimes_{r}H)\\ @A i_{*}' AA @|\\
KK_{*}(C_{0}(G\times_{H}X)\rtimes_{r} H, A \rtimes_{r}H)@>
i_{*}'([\lambda_{X}])\otimes\cdot
>> KK_{*}(\C, A \rtimes_{r}H)\\
@V p_{*}\circ q_{*} VV @V q_{*} VV\\
KK_{*}(C_{0}(G\times_{H}X)\rtimes_{r} G, A \rtimes_{r}
G)@>[\lambda_{G\times_{H}X}]\otimes\cdot
>> KK_{*}(\C, A \rtimes_{r}G)
\end{CD}\]


We finish the proof by splicing together these two diagrams.
\end{proof}
\begin{lem}
\label{**} Let $X$ be an $H$-compact subset of $P_c(H)$. Then the
following diagram commutes:
\[\begin{CD}
KK_{*}^{H}(C_{0}(X), A) @> \mu_{H} >> K_{*}(A \rtimes_{r}H)\\ @V
\mathcal{R}_{H,X}^G VV @V \mathcal{T}_{H}^{G} VV\\
KK_{*}^{G}(C_{0}(G\cdot X), A) @> \mu_{G} >> K_{*}(A
\rtimes_{r}G).
\end{CD}\]
\end{lem}
\begin{proof} By Lemma \ref{induced} we know that the induced space
$G \times_H X$ is $G$-homeomorphic to the $G$-saturation $G \cdot
X$. We now apply Lemma \ref{*}.
\end{proof}
\begin{thm}
\label{limittheorem} Let $G$ be a locally compact, second
countable Hausdorff topological group and let $A$ be a
$G-C^*$-algebra. Let $G$ be the union of open subgroups $G_{n}$
such that the Baum-Connes conjecture with coefficients $A$ is true
for each $G_{n}$. Then the Baum-Connes conjecture with
coefficients $A$ is true for $G$.
\end{thm}
\begin{proof}
\label{***} We start with the commutative diagram in Lemma
\ref{**} and take the direct limit over all $H$-compact subsets of
$P_c(H)$.   We then obtain the commutative diagram
\[
\xymatrix{
  \Kt{H,A} \ar[d]_{\mathcal{R}_{H}^{G}} \ar[r]^{\mu_{H}} & K_{*}(A \rtimes_{r} H)
  \ar[d]^{\mathcal{T}_{H}^{G}} \\
  \Kt{G,A} \ar[r]^{\mu_{G}} & K_{*}(A \rtimes_{r} G)}
\]

If $G$ is the union of open subgroups $G_{i}$ each of which
satisfies $BCC$, then applying Theorem \ref{Tmap} and Theorem
\ref{Rmap} along with the above commutative diagram is enough to
show that $G$ satisfies $BCC$.
\end{proof}

\begin{thm}
\label{result} Let $F$ be a global field, $\A$ its ring of adeles,
$G$ a linear algebraic group over $F$.  Let $F_v$ denote a place
of $F$. If BCC is true for each local group $G(F_v)$ then BCC is
true for the adelic group $G(\A)$.
\end{thm}

\begin{proof} Let $v_{1}, v_{2}, v_{3},\dots$ be an ordering of the
finite places of $F$, let $S_{\infty}$ denote the finite set of
all infinite places of $F$, and let $S(n) = \{v_1, v_2, \dots,
v_n\} \cup S_{\infty}$.  Let $G_n = G_{S(n)}$ in the notation of
section $2$.

If $\Gamma$ is a compact group then $\E \Gamma$ is a point, and we
have

\[\Kt{\Gamma, B} \cong KK^{\Gamma}_*(B) \cong K_*(B \rtimes \Gamma)\]

by the Green-Julg theorem.   Therefore BCC is true for any compact
group.   Now $G_n$ is a product of finitely many local groups and
one compact group.   But BCC is stable under the direct product of
finitely many groups, by an important result of Chabert-Echterhoff
\cite[Theorem 3.17]{CE}.   Therefore BCC is true for each {\it
open} subgroup $G_n$.  Now the adelic group $G(\A)$ is the
ascending union of the open subgroups $G_n$, therefore BCC is true
for $G(\A)$, by Theorem \ref{limittheorem}.
\end{proof}
\appendix
\section{Biller's theorem}
We give here the full statement of Theorem 3.8 in
Biller\cite{biller}.  The stabilizer of $x \in X$ is denoted
$G_x$.

{\bf Theorem} (Existence of slices). Let $G$ be a locally compact
group acting properly on a completely regular space $X$, and
choose $x\in X$. Then there is a convergent filter basis
$\mathcal{N}$ that consists of compact subgroups of $G$ normalized
by $G_x$ such that for every $N \in \mathcal{N}$, the coset space
$G/G_xN$ is a manifold and $x$ is contained in a $G_xN$-slice for
the action of $G$ on $X$.  In particular, some neighbourhood of
the orbit $G \cdot x$ is a locally trivial fibre bundle over the
manifold $G/G_xN$.

The dimension of $G \cdot x$ is infinite if and only if $N \in
\mathcal{N}$ may be chosen such that the dimension of $G/G_xN$ is
arbitrarily high. If the dimension of $G \cdot x$ is finite, then
$N \in \mathcal{N}$ may be chosen in such a way that $dim \,
G/G_xN = dim \, G \cdot x$.

Paul Baum, Mathematics Department, Pennsylvania State University,
University Park, PA 16802, USA. Email: baum@math.psu.edu

Stephen Millington, Compass Computer Consultants Ltd, Compass
House, Frodsham, Cheshire WA6 0AQ, UK. Email:
stem@quick.freeserve.co.uk

Roger Plymen, Mathematics Department, Manchester University,
Manchester M13 9PL, UK. Email: roger@ma.man.ac.uk
\end{document}